\documentclass[11pt]{amsart}

\usepackage{amsthm,amsmath,amssymb,amscd}
\usepackage{amssymb}
\usepackage{graphics}
\newtheorem{theorem}{Theorem}[section]
\newtheorem{corollary}[theorem]{Corollary}
\newtheorem{remark}[theorem]{Remark}
\newtheorem{lemma}[theorem]{Lemma}
\newtheorem{proposition}[theorem]{Proposition}
\newtheorem{definition}{Definition}[section]

\numberwithin{equation}{section}

\begin{document}

\title [Hamiltonian F-stability of complete Lagrangian self-shrinkers]
 {Hamiltonian F-stability of complete Lagrangian self-shrinkers}

\author{Liuqing Yang}

\address{Academy of Mathematics and Systems Sciences, Chinese Academy of Sciences, Beijing 100190, P. R. of China.}
\address{MIT, Dept. of Math., 77 Massachusetts Avenue, Cambridge, MA 02139-4307.}
\email{yangliuqing@amss.ac.cn}


\date{}

\maketitle

\begin{abstract}
In this paper, we study the Lagrangian F-stability and Hamiltonian F-stability of Lagrangian self-shrinkers. We prove a characterization theorem for the Hamiltonian F-stability of $n$-dimensional complete Lagrangian self-shrinkers without boundary, with polynomial volume growth and with the second fundamental form satisfying the condition that there exist constants $C_0>0$ and $\varepsilon<\frac{1}{16n}$ such that $|A|^2\leq C_0+\varepsilon |x|^2$. We characterize the Hamiltonian F-stablity by the eigenvalues and eigenspaces of the drifted Laplacian.
\end{abstract}

\vspace{.2in}


\section{Introduction}

\allowdisplaybreaks

\vspace{.1in}

\noindent

An $n$-dimensional submanifold $\Sigma^n$ of $\textbf{R}^{n+p}$ is called a self-shrinker if it is the time $t=-1$ slice of a self-shrinking mean curvature flow that disappears at $(0,0)$, i.e. of a mean curvature flow satisfying $\Sigma_t=\sqrt{-t}\Sigma_{-1}$. We can also consider a self-shrinker as a submanifold that satisfies
\begin{eqnarray*}
  H=-\frac{1}{2}x^\perp.
\end{eqnarray*}

Self-shrinkers are very important singularities of the mean curvature flow.

According to the blow up rate of the second
fundamental form, Huisken \cite{Hu2} classified the singularities of
mean curvature flows into two types: Type I and Type II. In 1984,
Huisken \cite{Hu1} showed that, if the initial hypersurface in
$\textbf{R}^{n+1}$ is strictly convex, then along the mean curvature
flow, the surface will be strictly convex at each time, and the mean
curvature flow will contract to a point at a finite time $T$.
Moreover, the normalized mean curvature flow will converge to a
round sphere. In 1990, Huisken \cite{Hu2} proved that any Type I
singularity of the mean curvature flow must be a self-shrinker by
using the monotonicity formula. He also proved that the only compact
self-shrinkers with nonnegative mean curvature are spheres.

In \cite{CM}, Colding-Minicozzi introduced the concept of F-stability and entropy-stability of a self-shrinker, and gave a classification of self-shrinkers in the hypersurface case. The definitions of many concepts in their paper can be naturally generalized to the higher codimension case (cf. \cite{ALW,AS,LL}).

Given $x_0\in \textbf{R}^{n+p}$ and $t_0>0$, $F_{x_0, t_0}$ is defined by
\begin{eqnarray*}
  F_{x_0,t_0}(\Sigma)=(4\pi t_0)^{-\frac{n}{2}}\int_\Sigma e^{-\frac{|x-x_0|^2}{4t_0}}d\mu.
\end{eqnarray*}
The entropy $\lambda=\lambda(\Sigma)$ is the supremum of the $F_{x_0, t_0}$ functionals
\begin{eqnarray*}
  \lambda=\sup_{x_0\in \textbf{R}^{n+p}, t_0>0} F_{x_0, t_0}(\Sigma).
\end{eqnarray*}
In \cite{CM}, Colding-Minicozzi proved that self-shrinkers are the critical points for the $F_{0,1}$ functional by computing the first variation formula of $F_{0,1}$. They also computed the second variation formula, and defined F-stability of a self-shrinker by modding out translations. More precisely, a self-shrinker $\Sigma$ is said to be F-stable if for every compactly supported variation $\Sigma_s$ with $\Sigma_0=\Sigma$, there exist variations $x_s$ of $x_0$ and $t_s$ of $t_0$ that make $F''=(F_{x_s, t_s}(\Sigma_s))''\geq 0$ at $s=0$. They also defined a self-shrinker to be entropy-stable if it is a local minimum for the entropy functional.

Colding-Minicozzi \cite{CM} showed that the round sphere and hyperplanes are the only F-stable self-shrinkers in $\textbf{R}^{n+1}$. By studying the relationship between F-stability and entropy-stability, they proved that every singularity other than spheres and cylinders can be perturbed away.

In 2002, Andrews-Li-Wei \cite{ALW}, Arezzo-Sun \cite{AS} and Lee-Lue \cite{LL} independently generalized Colding-Minicozzi's work \cite{CM} from the hypersurface case to the higher codimensional case. They computed the first and second variation formulae of the F-functional, and studied F-stability of self-shrinkers in higher codimension.

The normal bundle brings much difficulty for the classification of self-shrinkers in higher codimensions. In \cite{S}, Smoczyk classified self-shrinkers with $H\neq 0$ and parallel principal normal. Andrews-Li-Wei \cite{ALW} classified F-stable self-shrinkers with parallel principal normal by using Smoczyk's result. In \cite{LL}, Lee-Lue found an equivalent condition to the F-stabiliy. Moreover, they proved that in some cases the closed Lagrangian self-shrinkers given by Anciaux in \cite{A} are Lagrangian F-unstable. See \cite{ALW,AS,CL,CiL,LL,LW,S}, etc. for some other interesting results on self-shrinkers in higher codimensions.

A self-shrinker $\Sigma^n$ in $\textbf{C}^n$ is called a Lagrangian self-shrinker, if it is also a Lagrangian submanifold.
Lagrangian self-shrinkers are very important examples of self-shrinkers in higher codimension. Anciaux \cite{A}
and Joyce-Lee-Tsui \cite{JLT} constructed some examples of Lagrangian self-shrinkers.

In 1990, Oh \cite{O} introduced the notion of Hamiltonian stability of minimal Lagrangian submanifolds in K\"ahler-Einstein manifolds, which means stability under Hamiltonian variations. He found a criterion for the Hamiltonian stability, which reduces the Hamiltonian stability to the first eigenvalue of the Laplacian $\Delta$ acting on functions. Motivated by this result, we also get a characterization theorem for Hamiltonian F-stability (see Definition \ref{LHs}) of complete Lagrangian self-shrinkers in $\textbf{C}^n$, which reduces the Hamiltonian F-stability to the spectral property of drifted Laplacian $\mathcal{L}$ acting on functions. Since the F-stability is defined by modding out the translations, our result is related to not only eigenvalues but also eigenfunctions. Our theorem is as follows:
\begin{theorem}
   Let $\Sigma^n\subset \textbf{C}^n$ be a smooth complete Lagrangian self-shrinker without boundary and with polynomial volume growth. Suppose there exist constants $C_0>0$ and $\varepsilon<\frac{1}{16n}$ such that $|A|^2\leq C_0+\varepsilon |x|^2$. Then the following statements are equivalent:

  (i) $\Sigma$ is Hamiltonian F-stable.

  (ii) $\lambda_1({\mathcal{L}})=\frac{1}{2}$, $\lambda_2({\mathcal{L}})\geq 1$, and the eigenspace corresponding to the eigenvalue $\frac{1}{2}$ is spanned by coordinate functions.
\end{theorem}

\begin{remark}
  Many manifolds satisfy our condition on the second fundamental form that $|A|^2\leq C_0+\varepsilon|x|^2 (C_0>0, \varepsilon<\frac{1}{16n})$. For example, the clifford torus, cylinders, and manifolds that are asymptotic to cones, and so on.
\end{remark}

With the characterization theorem, it is natural to think about classifying Hamiltonian F-stable Lagrangian self-shrinkers. It is natural to ask whether Clifford torus is the only closed Hamiltonian F-stable Lagrangian self-shrinker in $\textbf{C}^n$, whether $\underbrace{S^1(\sqrt{2})\times\cdots S^1(\sqrt{2})}_k\times \textbf{R}^{n-k}$ ($0\leq k\leq n-1$) are the only noncompact Hamiltonian F-stable Lagrangian self-shrinkers without boundary, with polynomial volume growth and satisfying $|A|^2\leq C_0+\varepsilon|x|^2\ (C_0>0, \varepsilon<\frac{1}{16n})$ in $\textbf{C}^n$, and whether there exist any other examples of Hamiltonian F-stable Lagrangian self-shrinkers. In \cite{CM}, Colding-Minicozzi classified F-stable self-shrinkers by showing that F-stability implies mean convexity (i.e. $H\geq 0$) and then classifying the mean convex self-shrinkers. However, in our case, the method does not apply. In fact, it is hard to get information about the mean curvature from Hamiltonian F-stability. This makes classification very difficult.

After we got the characterization theorem for Hamiltonian F-stability of closed Lagrangian self-shrinkers, and while we were trying to prove some classification results, we found on December 18, 2013 on arXiv that Li-Zhang \cite{LZ}
also obtained the same characterization theorem for the closed case. For the purpose of completeness, and to make our paper more readable, we will also include the closed case and our own proof of it in this paper.

\vspace{.1in}

\noindent \emph{Acknowledgement: This work was completed while the author was a visiting PhD student at MIT. The author is grateful for the facilities provided. She would like to thank Professor Tobias Holck Colding for invitation and arranging this visit, and for his encouragement. She would like to thank Professor William P. Minicozzi II for introducing this problem to her and for many very helpful discussions, suggestions and comments. She would also like to thank Professor Jiayu Li for his encouragement. She acknowledges the China Scholarship Council (Grant No. 201304910263) for supporting her visit to MIT.}

\vspace{.2in}

\section{Preliminaries}
In this section, we recall some known results that were first proved by Colding-Minicozzi \cite{CM} for the hypersurface case, and then generalized by \cite{ALW,AS,LL} to the higher codimension case. The results will be used in the following sections.

Recall that
given $x_0\in \textbf{R}^{n+p}$ and $t_0>0$, $F_{x_0, t_0}$ is defined by
\begin{eqnarray*}
  F_{x_0,t_0}(\Sigma)=(4\pi t_0)^{-\frac{n}{2}}\int_\Sigma e^{-\frac{|x-x_0|^2}{4t_0}}d\mu.
\end{eqnarray*}
The first variation formula of $F_{x_0, t_0}$ is
\begin{lemma}
  Let $\Sigma_s\subset\textbf{R}^{n+p}$ be a variation of $\Sigma$ with normal variation vector field $V$. If $x_s$ and $t_s$ are variations of $x_0$ and $t_0$ with $x_0'=y$ and $t_0'=h$, then
  \begin{eqnarray*}
    \frac{\partial}{\partial s}(F_{x_s, t_s}(\Sigma_s))&=&(4\pi t_0)^{-\frac{n}{2}}\int_{\Sigma}\left\{\left\langle H+\frac{(x-x_0)^\perp}{2t_0}, V\right\rangle+h\left(\frac{|x-x_0|^2}{4t_0^2}-\frac{n}{2t_0}\right)\right.\\
    &&\left.\ \ \ \ \ \ \ \ \ \ \ \ \ \ \ \ \ +\frac{\langle x-x_0, y\rangle}{2t_0}\right\}e^{-\frac{|x-x_0|^2}{4t_0}}.
  \end{eqnarray*}
\end{lemma}
It follows that
\begin{proposition}
  $\Sigma$ is a critical point for $F_{x_0, t_0}$ if and only if $H=-\frac{(x-x_0)^\perp}{2t_0}$.
\end{proposition}
The second variation formula at a critical point is
\begin{theorem}
  Suppose that $\Sigma$ is complete, $\partial\Sigma=\emptyset$, $\Sigma$ has polynomial volume growth, and $\Sigma$ is a critical point for $F_{0,1}$. If $\Sigma_s$ is a normal variation of $\Sigma$, $x_s$, $t_s$ are variations with $x_0=0$ and $t_0=1$, and
\begin{eqnarray*}
  \partial_s\Big|_{s=0}\Sigma_s=V, \ \ \partial_s\Big|_{s=0}x_s=y, \ \ and \ \ \partial_s\Big|_{s=0}t_s=h,
\end{eqnarray*}
then setting $F''=\partial_{ss}\Big|_{s=0}(F_{x_s,t_s}(\Sigma_s))$ gives
\begin{eqnarray*}
  F''=(4\pi)^{-\frac{n}{2}}\int_\Sigma\left(-\langle V, LV\rangle+\langle V, y\rangle-h^2|H|^2-2h\langle H, V\rangle-\frac{1}{2}|y^\perp|^2\right)e^{-\frac{|x|^2}{4}}
\end{eqnarray*}
where
\begin{eqnarray*}
  LV&=&\Delta V-\frac{1}{2}\nabla_{x^T}^\perp V+\big\langle \langle A, V\rangle, A\big\rangle+\frac{1}{2}V\nonumber\\
  &=&\left(\Delta V^\alpha-\frac{1}{2}\langle x,\nabla V^\alpha\rangle+g^{ik}g^{jl}V^\beta h^\beta_{ij}h^\alpha_{kl}+\frac{1}{2}V^\alpha\right)e_\alpha.
\end{eqnarray*}
\end{theorem}

The linear operator defined by
\begin{eqnarray*}
  {\mathcal{L}}v=\Delta v-\frac{1}{2}\langle x, \nabla v\rangle=e^{\frac{|x|^2}{4}}div_\Sigma\left(e^{-\frac{|x|^2}{4}}\nabla v\right)
\end{eqnarray*}
is self-adjoint in a weighted $L^2$ space. This follows immediately from Stokes' theorem. More precisely,
\begin{lemma}
  If $\Sigma\subset \textbf{R}^{n+p}$ is a submanifold of $\textbf{R}^{n+p}$, $u$ is a $C^1$ function with compact support, and $v$ is a $C^2$ function, then
  \begin{eqnarray}
    \int_\Sigma u(\mathcal{L}v)e^{-\frac{|x|^2}{4}}=-\int_{\Sigma}\langle \nabla v, \nabla u\rangle e^{-\frac{|x|^2}{4}}
  \end{eqnarray}
\end{lemma}
\begin{corollary}\label{uLv}
  Suppose that $\Sigma\subset \textbf{R}^{n+p}$ is a complete submanifold of $\textbf{R}^{n+p}$ without boundary. If $u$, $v$ are $C^2$ functions with
  \begin{eqnarray*}
    \int_\Sigma(|u\nabla v|+|\nabla u||\nabla v|+|u{\mathcal{L}}v|)e^{-\frac{|x|^2}{4}}<\infty,
  \end{eqnarray*}
  then we get
  \begin{eqnarray*}
    \int_\Sigma u({\mathcal{L}}v)e^{-\frac{|x|^2}{4}}=-\int_\Sigma \langle \nabla v, \nabla u\rangle e^{-\frac{|x|^2}{4}}.
  \end{eqnarray*}
\end{corollary}

Now we display here some known properties that will be used in our paper.
Denote by $x^A$ ($A=1, 2, \cdots, x+p$) the coordinate functions of $\Sigma$ in $\textbf{R}^{n+p}$, i.e. $x^A$ is the $A$-th component of the position vector $x$, then
\begin{lemma} If $\Sigma^n\subset \textbf{R}^{n+p}$ is a self-shrinker, then
\begin{eqnarray}\label{LxA}
  {\mathcal{L}}x^A=-\frac{1}{2}x^A.
\end{eqnarray}
\begin{eqnarray}\label{deltax2}
  \Delta|x|^2=2n-|x^\perp|^2,
\end{eqnarray}
\begin{eqnarray}\label{Lx2}
  \mathcal{L}|x|^2=2n-|x|^2,
\end{eqnarray}
\begin{eqnarray}\label{LH}
  LH=H,
\end{eqnarray}
and for every constant vector field $y$,
\begin{eqnarray}\label{Ly}
  Ly^\perp=\frac{1}{2}y^\perp.
\end{eqnarray}
\end{lemma}
From (\ref{LxA}) we see the coordinate functions are eigenfunctions of $\mathcal{L}$ corresponding to the eigenvalue $\frac{1}{2}$. From (\ref{LH}) and (\ref{Ly}) we see $H$ and $y^\perp$ are both vector-valued eigenfunctions of $L$.

\section{Lagrangian F-stability and Hamiltonian F-stability}
In this section, we will define Lagrangian F-stability and Hamiltonian F-stability of Lagrangian self-shrinkers (also defined in \cite{AS,LL,LZ}).
Recall the definition of Lagrangian and Hamiltonian variations on a Lagrangian submanifold.
\begin{definition}\cite{O}
  Let $(M, \bar\omega)$ be a symplectic manifold $M$. Let $\Sigma\subset M$ be a Lagrangian submanifold and $V$ be a vector field along $\Sigma$. $V$ is called a Lagrangian (resp. Hamiltonian) variation if it satisfies that the one form $i^*(V\rfloor\bar\omega)$ on $\Sigma$ is closed (resp. exact).
\end{definition}
Note that both Lagrangian variations and Hamiltonian variations have some equivalent definitions.
\begin{lemma}\label{Lv}\cite{LL}
  A normal variation $V$ on $\Sigma$ is Lagrangian if and only if
  \begin{eqnarray*}
    \langle \nabla_X^\perp V, JY\rangle=\langle\nabla_Y^\perp V, JX\rangle,
  \end{eqnarray*}
  where $\nabla^\perp$ is the normal connection on $N\Sigma$ and $X, Y\in T\Sigma$.
\end{lemma}
\begin{lemma}\label{Hv}\cite{O}
  A normal variation $V$ on $\Sigma$ is Hamiltonian if and only if
  \begin{eqnarray*}
    V=J\nabla f,
  \end{eqnarray*}
  where $f$ is a function on $L$ and $\nabla$ is the gradient on $L$ with respect to the induced metric.
\end{lemma}

Now we are ready to define Lagrangian F-stability and Hamiltonian F-stability of Lagrangian self-shrinkers.
\begin{definition}\label{LHs}
  We say a Lagrangian self-shrinker $\Sigma$ is Lagrangian (resp. Hamiltonian) F-stable if for every compactly supported Lagrangian (resp. Hamiltonian)  variations $\Sigma_s$ with $\Sigma_0=\Sigma$, there exist variations $x_s$ of $x_0$ and $t_s$ of $t_0$ that make $F''\geq 0$.
\end{definition}

It is obvious that a Hamiltonian variation is also a Lagrangian variation, so Lagrangian F-stability implies Hamiltonian F-stability.

\section{Examples}
We begin with some simple examples to observe the Lagrangian F-stability and Hamiltonian F-stability of Lagrangian self-shrinkers.  It is well known that the simplest example of an $n$-dimensional closed Lagrangian self-shrinkers in $\textbf{C}^n$ is the Clifford torus  $T^n$, while the simplest example of an $n$-dimensional noncompact Lagrangian self-shrinker in $\textbf{C}^n$ besides $\textbf{R}^n$ is the cylinder $S^1(\sqrt{2})\times \textbf{R}^{n-1}$ . In this section, we will study the Lagrangian F-stability and Hamiltonian F-stability of these two examples. The proof in this section is inspired by Colding-Minicozzi's proof of Lemma 4.23 in \cite{CM}.

\begin{theorem}\label{cylinder}
  The cylinder $S^1(\sqrt{2})\times \textbf{R}^{n-1}\subset\textbf{C}^n$ is Lagrangian F-stable, hence it is also Hamiltonian F-stable.
\end{theorem}
{\it Proof. } For simplicity, we only prove the case $n=2$. The proof for the case $n>2$ is very similar.
The cylinder can be expressed as
\begin{eqnarray*}
  F(\theta, t)=\left(\sqrt{2}\cos\theta,\sqrt{2}\sin\theta, t, 0\right).
\end{eqnarray*}
We choose
\begin{eqnarray*}
  e_1=F_\theta=(-\sqrt{2}\sin\theta, \sqrt{2}\cos\theta, 0, 0), \ \ e_2=F_t=(0, 0, 1, 0),
\end{eqnarray*}
then
\begin{eqnarray*}
  (g_{ij})_{1\leq i, j\leq 2}=\left(
                                                                                                         \begin{array}{cc}
                                                                                                           2 & 0 \\
                                                                                                           0 & 1 \\
                                                                                                         \end{array}
                                                                                                       \right).
\end{eqnarray*}
Using the standard complex structure $J$ in $\textbf{C}^2$,
\begin{eqnarray*}
  J=\left(
      \begin{array}{cccc}
        0 & -1 & 0 & 0 \\
        1 & 0 & 0 & 0 \\
        0 & 0 & 0 & -1 \\
        0 & 0 & 1 & 0 \\
      \end{array}
    \right),
\end{eqnarray*}
we choose
\begin{eqnarray*}
  e_3=\frac{Je_1}{|Je_1|}=(-\cos\theta, -\sin\theta, 0, 0),\ \  e_4=\frac{Je_2}{|Je_2|}=(0, 0, 0, 1).
\end{eqnarray*}
Thus $\{e_1, e_2\}$ is an orthogonal basis of $T\Sigma$, $\{e_3,e_4\}$ is an orthonormal basis of $N\Sigma$.
Moreover,
\begin{eqnarray*}
  F_{\theta\theta}=(-\sqrt{2}\cos\theta, -\sqrt{2}\sin\theta, 0, 0),\ \ F_{\theta t}=0, \ \ F_{tt}=0.
\end{eqnarray*}
Thus
\begin{eqnarray*}
  h^3_{11}=\langle F_{\theta\theta}, e_3\rangle=\sqrt{2}, \ \ other\  h^{\alpha}_{ij}=0.
\end{eqnarray*}
Hence
\begin{eqnarray*}
  H=g^{11}h^3_{11}e_3=\frac{\sqrt{2}}{2}e_3, \ \ H=-\frac{1}{2}F^\perp,\ \  |H|=\frac{\sqrt{2}}{2}.
\end{eqnarray*}
Now suppose $V=fe_3+ge_4$. By Lemma \ref{Lv} we get $V$ is Lagrangian if and only if
\begin{eqnarray}\label{Lagrangian}
\sqrt{2}f_t=g_\theta.
\end{eqnarray}
It is easy to check that $\nabla^\perp e_3=\nabla^\perp e_4=0$, thus
\begin{eqnarray}\label{VLVc}
  -\langle V, LV\rangle&=&-V^3LV^3-V^4LV^4\nonumber\\
  &=&-V^3{\mathcal{L}V^3}-V^3(g^{11})^2V^3(h^3_{11})^2-\frac{1}{2}(V^3)^2-V^4{\mathcal{L}}V^4-\frac{1}{2}(V^4)\nonumber\\
  &=&-f{\mathcal{L}}f-\frac{1}{2}f^2-\frac{1}{2}f^2-g{\mathcal{L}}g-\frac{1}{2}g^2\nonumber\\
  &=&-f\mathcal{L}f-f^2-g\mathcal{L}g-\frac{1}{2}g^2.
\end{eqnarray}
Therefore,
\begin{eqnarray*}
  F''&=&(4\pi)^{-\frac{n}{2}}\int_{S^1(\sqrt{2})\times \textbf{R}^1}\left\{-f(\mathcal{L}f+f)-g\left(\mathcal{L}g+\frac{1}{2}g\right)+f\langle y, e_3\rangle+g\langle y, e_4\rangle-\frac{1}{2}h^2-\sqrt{2}hf\right.\\
  &&\left.\ \ \ \ \ \ \ \ \ \ \ \ \ \ \ \ \ \ \ \ \ \ \ \ \ \ -\frac{1}{2}\langle y, e_3\rangle^2-\frac{1}{2}\langle y, e_4\rangle^2\right\}e^{-\frac{t^2+2}{4}}.
\end{eqnarray*}
We compute
\begin{eqnarray*}
  \mathcal{L}=\frac{1}{2}\frac{\partial^2}{\partial\theta^2}+\frac{\partial^2}{\partial t^2}-\frac{1}{2}t\frac{\partial}{\partial t}=\Delta_{S^1(\sqrt{2})}+\mathcal{L}_{\textbf{R}^1}.
\end{eqnarray*}
It is known that
\begin{eqnarray*}
  \lambda_k\left(\Delta_{S^1(\sqrt{2})}\right)=\frac{k^2}{2},
\end{eqnarray*}
with the associated eigenspace spanned by $\{\cos k\theta,\ \sin k\theta\}$, while
\begin{eqnarray*}
  \lambda_k(\mathcal{L}_{\textbf{R}^1})=\frac{k}{2},
\end{eqnarray*}
with the associated eigenspace spanned by Hermite Polynomials
\begin{eqnarray*}
  \left\{H_k(t)=(-1)^k e^\frac{t^2}{4}\frac{d^k}{dt^k}e^{-\frac{t^2}{4}}\right\}.
\end{eqnarray*}
Therefore,
\begin{eqnarray*}
  \lambda_0(\mathcal{L})=0,
\end{eqnarray*}
with the associated eigenspace spanned by $\{1\}$;
\begin{eqnarray*}
  \lambda_1(\mathcal{L})=\frac{1}{2},
\end{eqnarray*}
with the associated eigenspace spanned by $\{\cos\theta, \sin\theta, t\}$; and
\begin{eqnarray*}
  \lambda_2(\mathcal{L})=1,
\end{eqnarray*}
with the associated eigenspace spanned by $\{t\cos\theta, t\sin\theta, t^2-2\}$.
Notice that
\begin{eqnarray*}
  \int_{S^1(\sqrt{2})\times\textbf{R}^1} f\cdot t e^{-\frac{|x|^2}{4}}&=&\int_{S^1(\sqrt{2})}\int_{\textbf{R}^1} f\cdot t e^{-\frac{t^2+2}{4}}dtd\sigma\\&=&\int_{S^1(\sqrt{2})}\left(\int_{\textbf{R}^1} f\cdot(-2\mathcal{L}_{\textbf{R}^1} t) e^{-\frac
  {t^2}{4}}dt\right)e^{-\frac{1}{2}}d\sigma\\
  &=&2\int_{S^1(\sqrt{2})}\left(\int_{\textbf{R}^1} f_t\cdot 1e^{-\frac{t^2}{4}}dt\right)e^{-\frac{1}{2}}d\sigma\\&=&\sqrt{2}\int_{S^1(\sqrt{2})}\int_{\textbf{R}^1}g_\theta e^{-\frac{t^2+2}{4}}dtd\sigma\\
  &=&\sqrt{2}\int_{\textbf{R}^1}\left(\int_{S^1(\sqrt{2})}g_\theta d\sigma\right)e^{-{\frac{t^2+2}{4}}}dt\\
  &=&0,
\end{eqnarray*}
where we used (\ref{Lagrangian}) in the fourth equality.

Therefore, we can choose $a_1$, $a_2$, $a_3\in \textbf{R}$, so that
\begin{eqnarray*}
  f=a_1+a_2\cos\theta+a_3\sin\theta+f_0\triangleq a_1+\langle z, e_3\rangle+f_0,
\end{eqnarray*}
where $z=(-a_2, -a_3)$, and $-f_0\mathcal{L}f_0\geq f_0^2$.
We can also choose $b_1\in \textbf{R}$ so that
\begin{eqnarray*}
  g=b_1+g_0,\ \ -g_0\mathcal{L}g_0\geq\frac{1}{2}g_0^2.
\end{eqnarray*}
It follows that
\begin{eqnarray*}
  F''&\geq&(4\pi)^{-\frac{n}{2}}\int_{S^1(\sqrt{2})\times \textbf{R}^1}\left\{-a_1^2-\frac{1}{2}\langle z, e_3\rangle^2-\frac{1}{2}b_1^2+\langle z, e_3\rangle\langle y, e_3\rangle+b_1\langle y, e_4\rangle -\frac{1}{2}h^2-\sqrt{2}ha_1\right.\\
  &&\left.\ \ \ \ \ \ \ \ \ \ \ \ \ \ \ \ \ \ \ \ \ \ \ \ \ \ -\frac{1}{2}\langle y, e_3\rangle^2-\frac{1}{2}\langle y, e_4\rangle^2\right\}e^{-\frac{t^2+2}{4}}\\
  &=&(4\pi)^{-\frac{n}{2}}\int_{S^1(\sqrt{2})\times \textbf{R}^1}\left\{-(a_1+\frac{\sqrt{2}}{2}h)^2-\frac{1}{2}(\langle z, e_3\rangle-\langle y, e_3\rangle)^2-\frac{1}{2}(b_1-\langle y, e_4\rangle)^2\right\}e^{-\frac{t^2+2}{4}}.
\end{eqnarray*}
Choose $h=-\sqrt{2}a_1$, $\langle y, e_3\rangle=\langle z, e_3\rangle$, $\langle y, e_4\rangle=b_1$, i.e.
$h=-\sqrt{2}a_1$, $y=(-a_2, -a_3, 0, b_1)$, then $F''\geq 0$. Therefore $S^1(\sqrt{2})\times \textbf{R}^1$ is Lagrangian F-stable, thus also Hamiltonian F-stable.

\hfill Q.E.D.

\begin{theorem}
  The Clifford torus  $T^n=\underbrace{S^1(\sqrt{2})\times\cdots\times S^1(\sqrt{2})}_n\subset\textbf{C}^n$ is Hamiltonian F-stable, but Lagrangian F-unstable.
\end{theorem}
{\it Proof. } For simplicity, we only prove the case $n=2$. The proof for the case $n>2$ is similar.

The Clifford torus can be expressed as
\begin{eqnarray*}
  F(\theta, \varphi)=(\sqrt{2}\cos\theta, \sqrt{2}\sin\theta, \sqrt{2}\cos\varphi, \sqrt{2}\sin\varphi).
\end{eqnarray*}
We choose
\begin{eqnarray*}
  e_1=F_\theta=(-\sqrt{2}\sin\theta, \sqrt{2}\cos\theta, 0, 0), \ \ e_2=F_\varphi=(0, 0, -\sqrt{2}\sin\varphi, \sqrt{2}\cos\varphi),
\end{eqnarray*}
\begin{eqnarray*}
  e_3=\frac{Je_1}{|Je_1|}=(-\cos\theta, -\sin\theta, 0, 0), \ \ e_4=\frac{Je_2}{|Je_2|}=(0, 0, -\cos\varphi, -\sin\varphi).
\end{eqnarray*}
It is easy to compute that
\begin{eqnarray*}
  H=\frac{\sqrt{2}}{2}e_3+\frac{\sqrt{2}}{2}e_4, |H|=1.
\end{eqnarray*}

(i) Now suppose $V$ is a Hamiltonian variation. Then by Lemma \ref{Hv}, there exists a function $f$ such that
\begin{eqnarray*}
  V=\sqrt{2}J\nabla f=f_\theta e_3+f_\varphi e_4.
\end{eqnarray*}
It is easy to check that $\nabla^\perp e_3=0$, $\nabla^\perp e_4=0$. By similar computations with (\ref{VLVc}), we have
\begin{eqnarray}\label{VLVt}
  -\langle V, LV\rangle=-f_\theta\mathcal{L}f_\theta-f_{\theta}^2-f_\varphi\mathcal{L}f_\varphi-f_\varphi^2,
\end{eqnarray}
and
\begin{eqnarray*}
  F''&=&(4\pi)^{-\frac{n}{2}}\int_\Sigma\left\{-f_\theta(\mathcal{L}f_\theta+f_\theta)-f_\varphi(\mathcal{L}f_\varphi+f_\varphi)+f_\theta\langle y, e_3\rangle+f_\varphi\langle y, e_4\rangle-h^2-\frac{1}{2}\langle y, e_3\rangle^2\right.\\
  &&\left.\ \ \ \ \ \ \ \ \ \ \ \ \ \ \ \ -\frac{1}{2}\langle y, e_4\rangle^2\right\}e^{-1}.
\end{eqnarray*}
By the same method as in the proof of Theorem \ref{cylinder}, we can get the Hamiltonian F-stability of the Clifford torus.

(ii)It is easy to check that $V=e_3-e_4$ is a Lagrangian variation, and $F''<0$ for every $h$ and $y$. Therefore, $S^1(\sqrt{2})\times S^1(\sqrt{2})$ is Lagrangian F-unstable.

We leave the details to the readers. \hfill Q.E.D.

The Hamiltonian F-stability of the Clifford torus and the cylinder is also an immediate corollary of our characterization theorem for Hamiltonian F-stability of complete Lagrangian self-shrinkers. See section 7 and section 9 for more details.

The Lagrangian F-stability of the cylinder was also mentioned by Li-Zhang in \cite{LZ}. They also proved the results that the Clifford torus is Hamiltonian F-stable and Lagrangian F-unstable as a corollary of their theorems. See \cite{LZ} for more discussions on Lagrangian F-stability and Hamiltonian F-stability of closed Lagrangian self-shrinkers.

\section{The operator $L$ preserves Hamiltonian}
Note that in the previous section, when we were computing $LV$ in both examples, we used $\nabla^\perp e_3=\nabla^\perp e_4=0$. Then the computation of $LV$ was reduced to the action of ${\mathcal{L}}$ on coefficient functions. However, we cannot always find such a good frame for general Lagrangian self-shrinkers. So it seems not easy to compute $LV$ in the general case.

In \cite{O}, Oh studied Hamiltonian stability of minimal Lagrangian submanifolds in K\"ahler-Einstein manifolds, and characterized Hamiltonian stability by a condition on the first eigenvalue of $\Delta$ acting on functions. The key point of Oh's proof is that, for a minimal Lagrangian submanifold of a K\"ahler-Einstein manifold, the set of Hamiltonian variations is an invariant subspace of the Jacobi operator. It is natural to think that this propery also holds for Lagrangian self-shrinkers. Oh proved this property by using the isomorphism between sections of $NL$ and one-forms. He transferred all the computations from normal vector fields to forms, and used Hodge-decomposition. In \cite{LZ}, Li-Zhang also made their computations on forms, and used the twisted Hodge Laplacian and twisted Hodge-decomposition. However, the property inspired us to to show the following equality on vector fields directly, though the computations are essentially equivalent with those on forms. The equality well characterizes how the operator $L$ acts on Hamiltonian variations.

\begin{theorem} Suppose $\Sigma^n\subset\textbf{C}^n$ is a Lagrangian self-shrinker. Then for every function $f$ on $\Sigma$,
  \begin{eqnarray}\label{e3.1}
    LJ\nabla f=J\nabla ({\mathcal{L}} f+ f).
  \end{eqnarray}
  This implies that the set of Hamiltonian variations is an invariant subspace of the operator $L$.
\end{theorem}

{\it Proof. } Fix a point $p$. We choose a local orthonormal basis $\{e_i\}_{i=1}^n$ of $T\Sigma$ such that $\nabla_{e_i}e_j(p)=0$. Then since $\Sigma$ is Lagrangian, $\{e_{n+i}=Je_i\}_{i=1}^n$ is a local orthonomal basis of $N\Sigma$. In the following we compute at the point $p$. It is easy to compute that
\begin{eqnarray*}
LJ\nabla f&=&\Delta^\perp(J\nabla f)-\frac{1}{2}\nabla^\perp_{x^T}(J\nabla f)+h^{n+k}_{il}f_kh_{il}^{n+j}Je_j+\frac{1}{2}f_jJe_j\nonumber\\
&=&\left(f_{jii}-\frac{1}{2}\langle x, e_k\rangle f_{jk}+f_kh^{n+l}_{ik}h^{n+l}_{ij}+\frac{1}{2}f_j\right)Je_j,
\end{eqnarray*}
where in the last equality we used the Lagrangian property $h_{il}^{n+k}=h_{ik}^{n+l}$.
On the other hand,
\begin{eqnarray*}
  J\nabla({\mathcal{L}}f+f)&=&J\nabla\left(\Delta f-\frac{1}{2}x^T f+f\right)\nonumber\\
  &=&f_{iij}Je_j-\frac{1}{2}J\nabla\langle x^T,\nabla f\rangle+f_jJe_j\nonumber\\
  &=&\left(f_{iji}-f_iR_{jkik}-\frac{1}{2}e_j\langle x^T, \nabla f\rangle+f_j\right)Je_j\nonumber\\
  &=&\left(f_{jii}-f_ih_{ij}^{n+l}h_{kk}^{n+l}+f_i h_{jk}^{n+l}h_{ik}^{n+l}-\frac{1}{2}\langle\nabla_{e_j}x^T, \nabla f\rangle-\frac{1}{2}\langle
  x^T, \nabla_{e_j}\nabla f\rangle+f_j\right)Je_j\nonumber\\
  &=&\left(f_{jii}+\frac{1}{2}f_ih_{ij}^{n+l}\langle x^\perp, e_{n+l}\rangle+f_kh_{ik}^{n+l}h_{ij}^{n+l}-\frac{1}{2}\langle\overline{\nabla}_{e_j}x,\nabla f\rangle+\frac{1}{2}\langle\overline{\nabla}_{e_j}x^\perp, \nabla f\rangle\right.\nonumber\\
  &&\left.-\frac{1}{2}f_{jk}\langle x, e_k\rangle+f_j\right) Je_j\nonumber\\
  &=&\left(f_{jii}+\frac{1}{2}f_ih_{ij}^{n+l}\langle x^\perp, e_{n+l}\rangle+f_kh_{ik}^{n+l}h_{ij}^{n+l}-\frac{1}{2}\langle e_j, \nabla f\rangle-\frac{1}{2}\left\langle x^\perp, \overline{\nabla}_{e_j}(f_ke_k)\right\rangle\right.\nonumber\\
  &&\left.-\frac{1}{2}\langle x, e_k\rangle f_{jk}+f_j\right) Je_j\nonumber\\
  &=&\left(f_{jii}+\frac{1}{2}f_ih_{ij}^{n+l}\langle x^\perp, e_{n+l}\rangle+f_kh_{ik}^{n+l}h_{ij}^{n+l}-\frac{1}{2}f_j-\frac{1}{2}f_k\langle x^\perp, h_{jk}^{n+l}e_{n+l}\rangle\right.\nonumber\\
  &&\left.-\frac{1}{2}\langle x, e_k\rangle f_{jk}+f_j\right)Je_j\nonumber\\
  &=&\left(f_{jii}-\frac{1}{2}\langle x, e_k\rangle f_{jk}+f_k h_{ik}^{n+l}h_{ij}^{n+l}+\frac{1}{2}f_j\right) Je_j,
\end{eqnarray*}
where in the third equality we used the Ricci formula; in the fourth equality we used the Gauss equation; and in the fifth equality we used the self-shrinker equation $H=-\frac{1}{2}x^\perp$.
This proves the theorem.  \hfill Q.E.D.

\section{Second variation formula under Hamiltonian variations}
From this section to the end of this paper, we will use square brackets $[\cdot]$ to denote weighted integrals
\begin{eqnarray}\label{notation}
  [f]=(4\pi)^{-\frac{n}{2}}\int_\Sigma f e^{-\frac{|x|^2}{4}}.
\end{eqnarray}
Then the second variation formula of $F$-functional can be written as
\begin{eqnarray}\label{e4.1}
  F''=\left[-\langle V, LV\rangle+\langle V, y\rangle-h^2|H|^2-2h\langle H, V\rangle-\frac{1}{2}|y^\perp|^2\right],
\end{eqnarray}
where
\begin{eqnarray*}
  L=\Delta^\perp-\frac{1}{2}\nabla^\perp_{x^T}+\big\langle\langle A, \cdot\rangle, A\big\rangle+\frac{1}{2}.
\end{eqnarray*}
Now we assume $V$ is a Hamiltonian variation, that is, there exists a function $f$, such that $V=J\nabla f$. Putting it into (\ref{e4.1}), we have
\begin{eqnarray}\label{e4.2}
  F''=\left[-\langle J\nabla f, LJ\nabla f\rangle+\langle J\nabla f, y\rangle-h^2|H|^2-2h\langle H, J\nabla f\rangle-\frac{1}{2}|y^\perp|^2\right]
\end{eqnarray}
First note that,
\begin{lemma}For every smooth function $f$ on a Lagrangian self-shrinker, we have
  \begin{eqnarray}\label{e4.3}
    [\langle H, J\nabla f\rangle]=0.
  \end{eqnarray}
\end{lemma}
{\it Proof.}
By choosing a geodesic frame, we compute
\begin{eqnarray*}
  div_\Sigma(Jx)^T&=&e_k\langle Jx, e_k\rangle=\langle\overline{\nabla}_{e_k}Jx, e_k\rangle+\langle Jx, \overline{\nabla}_{e_k}e_k\rangle=\langle J\overline{\nabla}_{e_k}x, e_k\rangle+\langle Jx, H\rangle\\
  &=&\langle J e_k, e_k\rangle+\langle Jx, H\rangle=\langle Jx, H\rangle,
\end{eqnarray*}
thus
\begin{eqnarray*}
  e^{-\frac{|x|^2}{4}}div_\Sigma\left(f(Jx)^Te^{-\frac{|x|^2}{4}}\right)&=&\langle \nabla f, (Jx)^T\rangle+f\langle Jx, H\rangle-\frac{1}{2}f\langle (Jx)^T, x^T\rangle\\
  &=&\langle \nabla f, J(x^\perp)\rangle-\frac{1}{2}f\langle Jx, x^\perp\rangle-\frac{1}{2}f\langle Jx, x^T\rangle\\
  &=&-2\langle\nabla f, JH\rangle-\frac{1}{2}f\langle Jx, x\rangle\\
  &=&2\langle J\nabla f, H\rangle.
\end{eqnarray*}
By approximation and Stokes' theorem, we get (\ref{e4.3}). \hfill Q.E.D.

Substituting (\ref{e3.1}) and (\ref{e4.3}) into (\ref{e4.2}), we immediately get the second variation formula under Hamiltonian variations.
\begin{theorem}
  Suppose that $\Sigma$ is complete, $\partial\Sigma=\emptyset$, $\Sigma$ has polynomial volume growth, and $\Sigma$ is a critical point for $F_{0,1}$. If $\Sigma_s$ is a Hamiltonian variation of $\Sigma$, $x_s$, $t_s$ are variations with $x_0=0$ and $t_0=1$, and
\begin{eqnarray*}
  \partial_s\big|_{s=0}\Sigma_s=V=J\nabla f, \ \ \partial_s\Big|_{s=0}x_s=y, \ \ and \ \ \partial_s\Big|_{s=0}t_s=h,
\end{eqnarray*}
Then setting $F''=\partial_{ss}\Big|_{s=0}(F_{x_s,t_s}(\Sigma_s))$ gives
  \begin{eqnarray}\label{F2H}
    F''=\left[-\langle \nabla f, \nabla ({\mathcal{L}}f+f)\rangle+\langle J\nabla f, y\rangle-h^2|H|^2-\frac{1}{2}|y^\perp|^2\right].
  \end{eqnarray}
\end{theorem}
Note that Li-Zhang (Proposition 4.1 in \cite{LZ}) also wrote the second variational formula under Hamiltonian variations, where they expressed it by forms and used $d_f^*$, the adjoint operator of $d$ in the weighted $L^2$ space.

\section{Hamiltonian F-stability of closed Lagrangian self-shrinkers}
In this section, we prove a characterization theorem for Hamiltonian F-stability of closed Lagrangian self-shrinkers. This theorem was also proved in \cite{LZ} (Theorem 1.3).
\begin{theorem}\label{closed}
  Suppose $\Sigma^n\subset \textbf{C}^n$ is a smooth closed Lagrangian self-shrinker, then the following statements are equivalent:

  (i) $\Sigma$ is Hamiltonian F-stable.

  (ii) $\lambda_1({\mathcal{L}})=\frac{1}{2}$, $\lambda_2({\mathcal{L}})\geq 1$, and the eigenspace corresponding to the eigenvalue $\frac{1}{2}$ is spanned by coordinate functions.
\end{theorem}
{\it Proof. }First we prove $(ii)\Rightarrow(i)$.
Given an arbitrary Hamiltonian vector field $V=J\nabla f$, by (ii) we can choose $a_0, a_A\in \textbf{R}$, $A=1,...,2n$, such that
\begin{eqnarray}\label{f}
  f=a_0+\sum_A a_Ax^A+f_0\triangleq a_0+\langle z, x\rangle+f_0,
\end{eqnarray}
where $z=(a_1, \cdots, a_{2n})$, $\mathcal{L}\langle z, x\rangle=-\frac{1}{2}\langle z, x\rangle$ and $-f_0{\mathcal{L}}f_0\geq f_0^2$.
It is easy to get
\begin{eqnarray}\label{e6.2}
  \nabla \langle z, x\rangle=z^T, \ \ J\nabla \langle z, x\rangle=Jz^T=(Jz)^\perp
\end{eqnarray}
Putting (\ref{f}) and (\ref{e6.2}) into (\ref{F2H}), and using the orthogonality of the different eigenspaces, we get that
\begin{eqnarray*}
  F''&\geq&\left[\frac{1}{2}|z^T|^2+\langle (Jz)^\perp, y^\perp\rangle-\frac{1}{2}|y^\perp|^2\right]\\
  &=&\left[-\frac{1}{2}\left|(Jz)^\perp-y^\perp\right|^2\right].
\end{eqnarray*}
We choose $h=0$ and $y=Jz$, then $F''\geq 0$. Therefore $\Sigma$ is Hamiltonian F-stable.

Now we prove $(i)\Rightarrow(ii)$. Assume the contrary that (ii) does not hold. Then either

(1) There exists a function $f$, such that ${\mathcal{L}}f=-\mu f$, where $0<\mu<1$, $\mu\neq\frac{1}{2}$; or

(2) There exists a function $f$, such that ${\mathcal{L}}f=-\mu f$, where $\mu=\frac{1}{2}$, $f$ is not a linear combination of coordinate functions, and \begin{eqnarray}\label{fxA}
  \left[f\cdot x^A\right]=0,\ \ A=1, \cdots, 2n.
\end{eqnarray}
First we consider case (1). Since $\mathcal{L}f=-\mu f$, by Corollary \ref{uLv}, we have
\begin{eqnarray*}
\left[|\nabla f|^2\right]=\left[-f\mathcal{L}f\right]=\mu\left[f^2\right].
\end{eqnarray*}
\begin{eqnarray*}
  LJ\nabla f=J\nabla(\mathcal{L}f+f)=J\nabla\left[(1-\mu)f\right]=(1-\mu)J\nabla f, \ \ 1-\mu\neq\frac{1}{2}.
\end{eqnarray*}
Since $Ly^\perp=\frac{1}{2}y^\perp$, and the eigenvector fields of $L$ corresponding to different eigenvalues are orthogonal with respect to the weighted $L^2$ inner product, we have
\begin{eqnarray*}
  \left[\langle J\nabla f, y\rangle\right]=\left[\langle J\nabla f, y^\perp\rangle\right]=0.
\end{eqnarray*}
Now we consider case (2). Then
\begin{eqnarray}\label{Lf}
\mathcal{L}f=-\frac{1}{2}f.
\end{eqnarray}
We compute that
\begin{eqnarray}\label{Jyx}
  J\nabla\langle Jy, x\rangle=J(Jy)^T=JJy^\perp=-y^\perp.
\end{eqnarray}
By (\ref{fxA}), we have
\begin{eqnarray}\label{fJyx}
  \left[f\langle Jy, x\rangle\right]=0
\end{eqnarray}
Hence by (\ref{Jyx}), (\ref{fJyx}) and using Corollary \ref{uLv}, we get
\begin{eqnarray*}
  \left[\langle J\nabla f, y\rangle\right]=\left[\langle J\nabla f, y^\perp\rangle\right]=\left[-\big\langle \nabla f, \nabla\langle Jy, x\rangle\big\rangle\right]=\left[{\mathcal{L}}f\cdot \langle Jy, x\rangle\right]=\left[-\frac{1}{2}f\langle Jy, x\rangle\right]=0
\end{eqnarray*}
Thus, in either case we have
\begin{eqnarray}\label{Jfy}
  \left[\langle J\nabla f, y\rangle\right]=0.
\end{eqnarray}
Therefore,
\begin{eqnarray*}
  F''&=&\left[\langle\nabla f, (\mu-1)\nabla f\rangle-h^2|H|^2-\frac{1}{2}|y^\perp|^2\right]=\left[(\mu-1)|\nabla f|^2-h^2|H|^2-\frac{1}{2}|y^\perp|^2\right]\\
  &=&\left[\mu(\mu-1)f^2-h^2|H|^2-\frac{1}{2}|y^\perp|^2\right]\leq\left[\mu(\mu-1)f^2\right]<0.
\end{eqnarray*}
This contradicts with (i).
This proves the theorem. \hfill Q.E.D.

Since the eigenvalues and eigenspaces of $\mathcal{L}$ on the Clifford torus satisfy (ii), it follows immediately that
\begin{corollary}
  The Clifford torus $T^n\subset\textbf{C}^n$ is Hamiltonian F-stable.
\end{corollary}

\section{Analysis on complete self-shrinkers}
In this section, we deduce some estimates on complete self-shrinkers, which will be used to prove our characterization theorem in the next section. The estimates we get in this section are motivated by section 3 in \cite{CM2}, where they got corresponding estimates on cylinders.

We use the same notations as in \cite{CM2}.
We denote the Gaussian $L^2$-norm
\begin{eqnarray}\label{L2}
\|u\|_{L^2}^2=\int_\Sigma u^2 e^{-\frac{|x|^2}{4}},
\end{eqnarray}
and the associated Gaussian $W^{1,2}$ and $W^{2,2}$ norms
\begin{eqnarray*}
  \|u\|_{W^{1,2}}^2=\int_\Sigma \left(u^2+|\nabla u|^2\right)e^{-\frac{|x|^2}{4}}
\end{eqnarray*}
and
\begin{eqnarray*}
  \|u\|_{W^{2,2}}^2=\int_\Sigma \left(u^2+|\nabla u|^2+|\nabla^2 u|^2\right)e^{-\frac{|x|^2}{4}}.
\end{eqnarray*}
\begin{lemma}\label{xu}
   Suppose $\Sigma^n\subset \textbf{C}^n$ is a smooth complete self-shrinker without boundary. Then there exists a constant $C$, such that if $u\in W^{1,2}(\Sigma)$, then
   \begin{eqnarray*}
     \||x|u\|_{L^2}^2\leq C\left(\|u\|_{L^2}^2+\|\nabla u\|_{L^2}^2\right)=C\|u\|_{W^{1,2}}^2.
   \end{eqnarray*}
\end{lemma}
{\it Proof. } It is easy to check that
\begin{eqnarray*}
  div_\Sigma x^T=n-\frac{1}{2}|x^\perp|^2,
\end{eqnarray*}
thus
\begin{eqnarray*}
  e^{\frac{|x|^2}{4}}div_{\Sigma}\left(u^2x^T e^{-\frac{|x|^2}{4}}\right)&=&2u\langle\nabla u, x^T\rangle+\left(n-\frac{1}{2}|x^\perp|^2\right)u^2-u^2\frac{|x^T|^2}{2}\nonumber\\
  &\leq&4|\nabla u|^2+\left(n-\frac{1}{2}|x^\perp|^2\right)u^2-u^2\frac{|x^T|^2}{4},
\end{eqnarray*}
where the inequality used the absorbing inequality $2ab\leq\frac{a^2}{4}+4b^2$.
By approximation, we can assume that $u$ has compact support on $\Sigma$, then by Stokes' theorem we have
\begin{eqnarray*}
  \frac{1}{4}\int_\Sigma u^2|x^T|^2e^{-\frac{|x|^2}{4}}\leq\int_\Sigma\left\{\left(n-\frac{1}{2}|x^\perp|^2\right)u^2+4|\nabla u|^2\right\}e^{-\frac{|x|^2}{4}}.
\end{eqnarray*}
The lemma follows since $|x|^2=|x^T|^2+|x^\perp|^2$. \hfill Q.E.D.

By induction, we have
\begin{lemma}\label{x2u}
     Suppose $\Sigma^n\subset \textbf{C}^n$ is a smooth complete self-shrinker without boundary. Then there exists a constant $C$, such that if $u\in W^{2,2}(\Sigma)$, then
   \begin{eqnarray}\label{e6.04}
     \||x|^2u\|_{L^2}^2\leq C\|u\|_{W^{2,2}}^2.
   \end{eqnarray}
\end{lemma}
{\it Proof. }By using (\ref{e6.04}) three times, we get
\begin{eqnarray*}
  \||x|^2u\|_{L^2}^2&=&\||x|(|x|u)\|_{L^2}^2\nonumber\\
  &\leq& C\||x|u\|_{W^{1,2}}^2\nonumber\\
  &=&C\left(\||x|u\|_{L^2}^2+\|\nabla(|x|u)\|_{L^2}^2\right)\nonumber\\
  &\leq& C\left(\|u\|_{W^{1,2}}^2+\|u\|_{L^2}^2+\||x||\nabla u|\|_{L^2}^2\right)\nonumber\\
  &\leq& C\left(\|u\|_{W^{1,2}}^2+\|\nabla u\|_{W^{1,2}}^2\right)\nonumber\\
  &\leq&C\|u\|_{W^{2,2}}^2,
\end{eqnarray*}
where the three inequalities used Lemma \ref{xu} for $|x|u$, $u$ and $|\nabla u|$, respectively.
This proves the lemma. \hfill Q.E.D.
\begin{lemma}\label{l6.3}
       Let $\Sigma^n\subset \textbf{C}^n$ be a smooth complete self-shrinker without boundary. Suppose there exist constants $C_0>0$ and $\varepsilon<\frac{1}{16n}$ such that $|A|^2\leq C_0+\varepsilon |x|^2$. If $u\in W^{1,2}(\Sigma)\cap C^2(\Sigma)$, and $\mathcal{L}u\in L^2(\Sigma)$, then $u\in W^{2,2}(\Sigma)$, and there exists a constant $C$, such that
   \begin{eqnarray}\label{e6.4}
     \|u\|_{W^{2,2}}^2\leq C\left(\|u\|_{L^2}^2+\|{\mathcal{L}}u\|_{L^2}^2\right).
   \end{eqnarray}
\end{lemma}
{\it Proof. } By integrating by parts, we get
\begin{eqnarray}\label{nablau}
  \|\nabla u\|_{L^2}^2=\left|\langle u, \mathcal{L}u\rangle_{L^2}\right|\leq \|u\|_{L^2}\|{\mathcal{L}}u\|_{L^2}\leq\frac{1}{2}\|u\|_{L^2}^2+\frac{1}{2}\|\mathcal{L}u\|_{L^2}^2.
\end{eqnarray}
It remains to bound $\|\nabla^2 u\|_{L^2}$. Let $\phi$ be a smooth function satisfying $\phi=1$ on $B_R$, $|\nabla \phi|\leq 1$ on $B_{R+2}\setminus B_R$, and $\phi=0$ on $\Sigma\setminus B_{R+2}$.
By direct computation, we get
\begin{eqnarray}\label{div}
  e^{\frac{|x|^2}{4}}div_\Sigma\left(\phi^2\{u_{ij}u_i-(\mathcal{L}u)u_j\}e^{-\frac{|x|^2}{4}}\right)&=&2\phi\phi_j\{u_{ij}u_i-(\mathcal{L}u)u_j
  \}\nonumber\\
  &&+\phi^2\left\{\frac{1}{2}\mathcal{L}|\nabla u|^2-(\mathcal{L}u)^2-\langle\nabla\mathcal{L}u, \nabla u\rangle\right\}.
\end{eqnarray}
We estimate
\begin{eqnarray}\label{phiphij}
  \left|\phi\phi_j\{u_{ij}u_i-(\mathcal{L}u)u_j\}\right|&\leq &|\phi||\nabla\phi||\nabla^2u||\nabla u|+|\phi||\nabla\phi||\mathcal{L}u||\nabla u|\nonumber\\
  &\leq&\delta\phi^2|\nabla^2u|^2+C_\delta|\nabla\phi|^2|\nabla u|^2+C\phi^2|\mathcal{L} u|^2,
\end{eqnarray}
where $\delta$ is to be determined later.

Recall the Bochner formula for the drifted Laplacian $\Delta_f u=\Delta u-\langle\nabla f, \nabla u\rangle$,
\begin{eqnarray*}
  \frac{1}{2}\Delta_f|\nabla u|^2=|\nabla^2 u|^2+\langle\nabla\Delta_f u, \nabla u\rangle +Ric_f (\nabla u, \nabla u).
\end{eqnarray*}
where $Ric_f=Ric+\nabla^2 f$ is the Bakry-\'Emery Ricci curvature. When $f=\frac{|x|^2}{4}$, $\Delta_f$ is just $\mathcal{L}$.
Then by our condition on the second fundamental form, we have
\begin{eqnarray}\label{Ricf}
  Ric_f(\nabla u, \nabla u)&=&\frac{1}{2}|\nabla u|^2-h^\alpha_{ik}h^\alpha_{jk}u_iu_j=\frac{1}{2}|\nabla u|^2-\left\{\sum_i\left(\Big(\sum_k h^\alpha_{ik}\Big)u_i\right)\right\}^2\nonumber\\
  &\geq&\frac{1}{2}|\nabla u|^2-\sum_i(\sum_k h^\alpha_{ik})^2\sum_i u_i^2\geq\frac{1}{2}|\nabla u|^2- n \sum_{i,k}(h^\alpha_{ik})^2\sum_i u_i^2\nonumber\\
  &=&\frac{1}{2}|\nabla u|^2-n|A|^2|\nabla u|^2 \nonumber\\
  &\geq& -C|\nabla u|^2-n\varepsilon |x|^2|\nabla u|^2.
\end{eqnarray}
Therefore, by the Bochner formula, we have
\begin{eqnarray}\label{BochnerL}
  \frac{1}{2}{\mathcal{L}}|\nabla u|^2-\langle \nabla{\mathcal{L}}u, \nabla u\rangle&=&|\nabla^2 u|^2+Ric_f(\nabla u, \nabla u)\nonumber\\
  &\geq& |\nabla^2 u|^2-C|\nabla u|^2-n\varepsilon|x|^2|\nabla u|^2.
\end{eqnarray}
Putting (\ref{phiphij}) and (\ref{BochnerL}) into (\ref{div}), we have
\begin{eqnarray}\label{e6.7}
  e^{\frac{|x|^2}{4}}div_\Sigma\left(\phi^2\{u_{ij}u_i-(\mathcal{L}u)u_j\}e^{-\frac{|x|^2}{4}}\right)&\geq&(1-\delta)\phi^2|\nabla^2 u|^2-C\phi^2|\nabla u|^2-C_\delta |\nabla\phi|^2|\nabla u|^2\nonumber\\
  &&-C\phi^2(\mathcal{L}u)^2-n\varepsilon\phi^2|x|^2|\nabla u|^2.
\end{eqnarray}
Thus by Stokes' theorem, we get
\begin{eqnarray}\label{deltaphi2}
  (1-\delta)\int_\Sigma\phi^2|\nabla^2 u|^2 e^{-\frac{|x|^2}{4}}&\leq& \int_\Sigma \left\{C\phi^2|\nabla u|^2+C_\delta|\nabla\phi|^2|\nabla u|^2+C\phi^2(\mathcal{L}u)^2\right.\nonumber\\&&\left.\ \ \ \ \ +n\varepsilon\phi^2|x|^2|\nabla u|^2\right\}e^{-\frac{|x|^2}{4}}.
\end{eqnarray}
Now we estimate $\int_\Sigma\phi^2|x|^2|\nabla u|^2 e^{-\frac{|x|^2}{4}}$.
Similar to the proof of Lemma \ref{xu}, we have
\begin{eqnarray*}
  e^{\frac{|x|^2}{4}}div_\Sigma\left(\phi^2|\nabla u|^2x^T e^{-\frac{|x|^2}{4}}\right)&=&2\phi\langle\nabla\phi, x^T\rangle|\nabla u|^2+2\phi^2|\nabla u|\langle\nabla|\nabla u|, x^T\rangle\\
  &&+\phi^2\left(n-\frac{1}{2}|x^\perp|^2\right)|\nabla u|^2-\phi^2|\nabla u|^2\frac{|x^T|^2}{2}\\
  &\leq&\frac{\delta}{4}\phi^2|\nabla u|^2|x^T|^2+C_\delta|\nabla \phi|^2|\nabla u|^2+4\phi^2|\nabla^2u|^2\\&&+\phi^2\left(n-\frac{1}{2}|x^\perp|^2\right)|\nabla u|^2-\phi^2|\nabla u|^2\frac{|x^T|^2}{4}.
\end{eqnarray*}
Thus by Stokes' Theorem, we get
\begin{eqnarray*}
  \frac{1-\delta}{4}\int_\Sigma\phi^2|x^T|^2|\nabla u|^2e^{-\frac{|x|^2}{4}}\leq\int_\Sigma \left\{C_\delta|\nabla\phi|^2|\nabla u|^2+\phi^2\left(n-\frac{1}{2}|x^\perp|^2\right)|\nabla u|^2+4\phi^2|\nabla^2u|^2\right\}e^{-\frac{|x|^2}{4}}.
\end{eqnarray*}
Therefore, by $|x|^2=|x^T|^2+|x^\perp|^2$, we have
\begin{eqnarray}\label{phi2}
  \ \ \ \ \ \ \  \int_\Sigma \phi^2|x|^2|\nabla u|^2e^{-\frac{|x|^2}{4}}\leq\int_\Sigma \left\{C_\delta|\nabla\phi|^2|\nabla u|^2+C_\delta\phi^2|\nabla u|^2+\frac{16}{1-\delta}\phi^2|\nabla^2 u|^2\right\}e^{-\frac{|x|^2}{4}}.
\end{eqnarray}
Combining (\ref{deltaphi2}) and (\ref{phi2}) gives
\begin{eqnarray}\label{nabla2u}
  &&\left(1-\delta-\frac{16n\varepsilon}{1-\delta}\right)\int_\Sigma\phi^2|\nabla^2 u|^2e^{-\frac{|x|^2}{4}}\nonumber\\&\leq&\int_\Sigma \left\{C_{\varepsilon, \delta}\phi^2|\nabla u|^2+C_{\varepsilon, \delta}|\nabla\phi|^2|\nabla u|^2+C\phi^2(\mathcal{L}u)^2\right\}e^{-\frac{|x|^2}{4}}.
\end{eqnarray}
Since $\varepsilon<\frac{1}{16n}$, we can choose $\delta>0$ so that $(1-\delta)^2>16n\varepsilon$.
Thus by (\ref{nablau}) and (\ref{nabla2u}), we get
\begin{eqnarray}\label{e6.8}
    \|\phi^2\nabla^2 u\|_{L^2}^2\leq C\left(\|u\|_{L^2}^2+\|\mathcal{L}u\|_{L^2}^2\right).
\end{eqnarray}
Thus the monotone convergence theorem gives $u\in W^{2,2}(\Sigma)$, and
\begin{eqnarray}
  \|\nabla^2u\|_{L^2}^2\leq C\left(\|u\|_{L^2}^2+\|\mathcal{L}u\|_{L^2}^2\right).
\end{eqnarray}
This proves the lemma. \hfill Q.E.D.

\section{Hamiltonian F-stability of complete Lagrangian self-shrinkers}
In this section, we prove our characterization theorem for Hamiltonian F-stability of complete Lagrangian self-shrinkers without boundary, with polynomial volume growth and with the second fundamental form satisfying the condition that there exist constants $C_0>0$ and $\varepsilon<\frac{1}{16n}$ such that $|A|^2\leq C_0+\varepsilon |x|^2$.

In \cite{CZ}, Cheng-Zhou studied the eigenvalues of the drifted Laplacian on complete metric measure spaces. In particular, they studied the spectrum of $\mathcal{L}$ on self-shrinkers, and proved that the spectrum of $\mathcal {L}$ is discrete for a properly immersed self-shrinker. Together with the result that for a self-shrinker, proper immersion, Euclidean volume growth, polynomial volume growth and finite weighted volume are equivalent each other (cf. \cite{CZ,DX}), their theorem can be stated as follows,
\begin{theorem}\cite{CZ}
  Let $\Sigma^n$ be a complete $n$-dimensional self-shrinker in the Euclidean space $\textbf{R}^{n+p}$, $p\geq 1$. Assume $\Sigma$ has polynomial volume growth, then the spectrum of $\mathcal{L}$ is discrete and consequently the first nonzero eigenvalue $\lambda_1$ of $\mathcal{L}$ has finite multiplicity and satisfies $\lambda_1\leq\frac{1}{2}$.
\end{theorem}

Since $\Sigma$ has polynomial volume growth, $0$ is the least eigenvalue of $\mathcal{L}$ with multiplicity one and the associated eigenfunctions are non-zero constant functions. Thus the set of all eigenvalues of $\mathcal{L}$ is an increasing sequence
\begin{eqnarray*}
  0=\lambda_0(\mathcal{L})<\lambda_1(\mathcal{L})<\lambda_2(\mathcal{L})<\cdots
\end{eqnarray*}
with $\lambda_i(\mathcal{L})\to\infty$ as $i\to\infty$. Moreover, by Lemma \ref{l6.3}, for each $i$, there exists a countable orthonormal base $\left\{\psi_j^{(i)}\right\}$ of $L^2(\Sigma)$ so that each $\psi_j^{(i)}\in W^{2,2}(\Sigma)$ is an eigenfunction of $\mathcal{L}$ associated with the eigenvalue $\lambda_i(\mathcal{L})$.
Our characterization theorem is as follows.
\begin{theorem}\label{complete}
  Let $\Sigma^n\subset \textbf{C}^n$ be a smooth complete Lagrangian self-shrinker without boundary and with polynomial volume growth. Suppose there exist constants $C_0>0$ and $\varepsilon<\frac{1}{16n}$ such that $|A|^2\leq C_0+\varepsilon |x|^2$. Then the following statements are equivalent:

  (i) $\Sigma$ is Hamiltonian F-stable.

  (ii) $\lambda_1({\mathcal{L}})=\frac{1}{2}$, $\lambda_2({\mathcal{L}})\geq 1$, and the eigenspace corresponding to the eigenvalue $\frac{1}{2}$ is spanned by coordinate functions.
\end{theorem}
{\it Proof. }First we prove $(ii)\Rightarrow(i)$.
Given an arbitrary compactly supported Hamiltonian vector field $V=J\nabla f$, by (ii) we can choose $a_0, a_A\in \textbf{R}$, $A=1,...,2n$, such that
\begin{eqnarray}\label{e6.1}
  f=a_0+\sum_{A=1}^{2n} a_Ax^A+f_0\triangleq a_0+\langle z, x\rangle+f_0,
\end{eqnarray}
where $z=(a_1, \cdots, a_{2n})$, $\mathcal{L}\langle z, x\rangle=-\frac{1}{2}\langle z, x\rangle$ and $-f_0{\mathcal{L}}f_0\geq f_0^2$.
Now the remaining part of the proof is essentially the same as the proof of $(ii)\Rightarrow(i)$ in Theorem \ref{closed}.

Now we prove $(i)\Rightarrow(ii)$. Since the eigenfunctions of $\mathcal{L}$ do not necessarily have compact support, we need to choose cutoff functions.
Let $\eta$ be a nonnegative smooth function on $[0, +\infty)$ satisfying
\begin{eqnarray*}
  \eta(s)=\left\{
            \begin{array}{ll}
              1, & \hbox{if $s\in [0, 1)$;} \\
              0, & \hbox{if $s\in [2, +\infty)$,}
            \end{array}
          \right.
\end{eqnarray*}
$|\eta'|\leq 2$, $|\eta''|\leq C$, and $|\eta'''|\leq C$.
Define a sequence of functions
\begin{eqnarray}\label{phij}
  \phi_j(x)=\eta\left(\frac{|x|^2}{j}\right).
\end{eqnarray}
Then $\phi_j\to 1$, and $|\nabla \phi_j|\to 0$ pointwise. Moreover,
\begin{eqnarray}\label{nablaphij}
  \nabla\phi_j=\eta'\frac{\nabla|x|^2}{j}=\eta'\frac{2x^T}{j}.
\end{eqnarray}
Using (\ref{deltax2}), we have
\begin{eqnarray}\label{deltaphij}
  \Delta\phi_j=\eta''\frac{\left|\nabla|x|^2\right|^2}{j^2}+\eta'\frac{\Delta|x|^2}{j}=\eta''\frac{4|x^T|^2}{j^2}+\eta'\frac{2n-|x^\perp|^2}{j}.
\end{eqnarray}
By (\ref{Lx2}), we get
\begin{eqnarray}\label{Lphij}
  \mathcal{L}\phi_j=\Delta\phi_j-\frac{1}{2}\langle\nabla\phi_j, x\rangle=\eta''\frac{4|x^T|^2}{j^2}+\eta'\frac{\mathcal{L}|x|^2}{j}=\eta''\frac{4|x^T|^2}{j^2}+\eta'\frac{2n-|x|^2}{j}.
\end{eqnarray}

Let $\{e_1, \cdots , e_n\}$ be a local orthonormal basis of $T\Sigma$, $\{e_{n+1}, \cdots, e_{2n}\}$ be a local
orthnormal basis of $N\Sigma$.
We make the following convention on the range of indices:
$1\leq i, k\leq n$; $n+1\leq\alpha \leq 2n.$
Recall that (cf. \cite{CL})
\begin{eqnarray*}
  x_{,i}=e_i,
\end{eqnarray*}
\begin{eqnarray*}
  x_{,ik}=h_{ik}^\alpha e_\alpha,
\end{eqnarray*}
and
\begin{eqnarray*}
  H^\alpha_{,i}=\frac{1}{2}h^\alpha_{ik}\langle x, e_k\rangle.
\end{eqnarray*}
It follows that
\begin{eqnarray*}
  |x|^2_{,ik}=2\langle x_{,i}, x_{,k}\rangle+2\langle x, x_{,ik}\rangle=2\langle e_i, e_k\rangle+2\langle x, h_{ik}^\alpha e_\alpha\rangle=2\delta_{ik}-4H^\alpha h_{ik}^\alpha,
\end{eqnarray*}
and
\begin{eqnarray*}
  \nabla|x^T|^2&=&\nabla|x|^2-\nabla|x^\perp|^2=2x^T-4\nabla|H|^2=2x^T-8H^\alpha H^\alpha_{,i}e_i\nonumber\\
&=&2x^T-4H^\alpha h^\alpha_{ik}\langle x, e_k\rangle e_i.
\end{eqnarray*}
Therefore,
\begin{eqnarray}\label{phijik}
  (\phi_j)_{ik}=\eta''\frac{4\langle{x, e_i}\rangle\langle{x, e_k}\rangle}{j^2}+\eta'\frac{|x|^2_{,ik}}{j}=\eta''\frac{4\langle{x, e_i}\rangle\langle{x, e_k}\rangle}{j^2}+\eta'\frac{2\delta_{ik}-4H^\alpha h_{ik}^\alpha}{j},
\end{eqnarray}
\begin{eqnarray}\label{nablaLphij}
  \nabla\mathcal{L}\phi_j&=&\eta'''\frac{4|x^T|^2}{j^2}\cdot\frac{\nabla|x|^2}{j}
+\eta''\frac{4\nabla|x^T|^2}{j^2}+\eta''\frac{2n-|x|^2}{j}\cdot\frac{\nabla|x|^2}{j}-\eta'\frac{\nabla|x|^2}{j}\nonumber\\
&=&\frac{8\eta'''|x^T|^2}{j^3}x^T+\eta''\frac{8x^T-16H^\alpha h^\alpha_{ik}\langle x, e_k\rangle e_i}{j^2}+\frac{2\eta''(2n-|x|^2)}{j^2}x^T-\frac{2\eta'}{j}x^T.
\end{eqnarray}
By our computations (\ref{phij})-(\ref{nablaLphij}), using the condition that $|A|^2\leq C_0+\varepsilon|x|^2$, we get
\begin{eqnarray}\label{phij2}
  |\phi_j|\leq 1,
\end{eqnarray}
\begin{eqnarray}\label{nablaphij2}
  |\nabla\phi_j|\leq C|x|,
\end{eqnarray}
\begin{eqnarray}\label{deltaphij2}
  |\Delta\phi_j|\leq C\left(1+|x|^2\right),
\end{eqnarray}
\begin{eqnarray}\label{Lphij2}
  |\mathcal{L}\phi_j|\leq C\left(1+|x|^2\right),
\end{eqnarray}
\begin{eqnarray}\label{hessphij}
  |\nabla^2\phi_j|\leq C\left(1+|x|^2\right),
\end{eqnarray}
and
\begin{eqnarray}\label{nablaLphij2}
  |\nabla\mathcal{L}\phi_j|\leq C\left(|x|+|x|^3\right).
\end{eqnarray}
 Now assume the contrary that (ii) does not hold. Then either

(1) There exists a function $f$, such that ${\mathcal{L}}f=-\mu f$, where $0<\mu<1$, $\mu\neq\frac{1}{2}$; or

(2) There exists a function $f$, such that ${\mathcal{L}}f=-\mu f$, where $\mu=\frac{1}{2}$, $f$ is not a linear combination of coordinate functions, and \begin{eqnarray}\label{fxA2}
  \left[f\cdot x^A\right]=0,\ \ A=1, \cdots, 2n.
\end{eqnarray}

Denote $f_j=\phi_j f$. In the following, we will use $V_j=J\nabla f_j$ as a variation. Then
\begin{eqnarray}\label{Fj}
  F_j''&=&\left[-\langle \nabla f_j, \nabla\left({\mathcal{L}}f_j+f_j\right)\rangle+\langle J\nabla f_j, y\rangle-h^2|H|^2-\frac{1}{2}|y^\perp|^2\right].
\end{eqnarray}
Direct computation gives
\begin{eqnarray}\label{fj1}
  \left[-\langle \nabla f_j, \nabla\left({\mathcal{L}}f_j+f_j\right)\rangle\right]&=&\left[-\langle \nabla (\phi_j f), \nabla\left({\mathcal{L}}(\phi_j f)+\phi_j f\right)\rangle\right]\nonumber\\
  &=&-\left[(1-\mu)|\phi_j|^2|\nabla f|^2+(2-2\mu)f\phi_j\langle\nabla f, \nabla\phi_j\rangle\right.\nonumber\\&&\left.+(1-\mu)f^2|\nabla\phi_j|^2+\phi_j\mathcal{L}\phi_j|\nabla f|^2+f\phi_j\langle\nabla f, \nabla\mathcal{L}\phi_j\rangle\right.\nonumber\\
&&\left.+f\mathcal{L}\phi_j\langle\nabla f, \nabla \phi_j\rangle+f^2\langle\nabla \phi_j, \nabla\mathcal{L}\phi_j\rangle\right.\nonumber\\
&&\left.+2\phi_j\langle\nabla f, \nabla\langle\nabla\phi_j, \nabla f\rangle\rangle+2f\langle\nabla\phi_j, \nabla\langle\nabla\phi_j, \nabla f\rangle\rangle\right].
\end{eqnarray}
From the proof of $(i)\Rightarrow(ii)$ in Theorem \ref{closed} we know
\begin{eqnarray*}
  \left[\langle J\nabla f, y^\perp\rangle\right]=0.
\end{eqnarray*}
Therefore
\begin{eqnarray}\label{fj2}
  \left[\langle J\nabla f_j, y\rangle\right]&=&\left[\langle J\nabla(\phi_jf), y^\perp\rangle\right]=\left[f\langle J\nabla\phi_j, y^\perp\rangle+\phi_j\langle J\nabla f, y^\perp\rangle\right]\nonumber\\
  &=&\left[f\langle J\nabla\phi_j, y^\perp\rangle+(\phi_j-1)\langle J\nabla f, y^\perp\rangle\right].
\end{eqnarray}
In order to use the dominated convergence theorem, we need to control all the terms that include $\phi_j$ in (\ref{fj1}) and (\ref{fj2}).

Since $f$ is an eigenfunction of $\mathcal{L}$, we know that $f\in W^{2,2}(\Sigma)$.
Note that our notations (\ref{notation}) and (\ref{L2}) satisfy
\begin{eqnarray*}
  \left[f^2\right]=(4\pi)^{-\frac{n}{2}}\|f\|_{L^2}^2.
\end{eqnarray*}
By Lemma \ref{xu} and Lemma \ref{x2u}, we have
\begin{eqnarray}\label{f1}
\left[|x|^2f^2\right]\leq C\|f\|_{W^{1,2}}^2,
\end{eqnarray}
and
\begin{eqnarray}
\left[|x|^4f^2\right]\leq C\|f\|_{W^{2,2}}^2.
\end{eqnarray}
Applying Lemma \ref{xu} on $u=|\nabla f|$ yields
\begin{eqnarray}\label{f2}
  \left[|x|^2|\nabla f|^2\right]\leq C\|\nabla f\|_{W^{1,2}}^2\leq C\|f\|_{W^{2,2}}^2.
\end{eqnarray}
Combining (\ref{f1})-(\ref{f2}) and the estimates (\ref{phij2})-(\ref{nablaLphij2}), we get
\begin{eqnarray}\label{Fj1}
  \left[|\phi_j|^2|\nabla f|^2\right]\leq\left[|\nabla f|^2\right]\leq C\|f\|_{W^{1,2}}^2,
\end{eqnarray}
\begin{eqnarray}
  \left[\left|f\phi_j\langle\nabla f, \nabla\phi_j\rangle\right|\right]\leq C\left[f|x||\nabla f|\right]\leq C\left[|x|^2f^2\right]+C\left[|\nabla f|^2\right]\leq C\|f\|_{W^{1,2}}^2,
\end{eqnarray}
\begin{eqnarray}
  \left[f^2|\nabla\phi_j|^2\right]\leq C\left[|x|^2f^2\right]\leq C\|f\|_{W^{1,2}}^2,
\end{eqnarray}
\begin{eqnarray}
  \left[\left|\phi_j\mathcal{L}\phi_j|\nabla f|^2\right|\right]\leq C\left[|\nabla f|^2\right]+C\left[|x|^2|\nabla f|^2\right]\leq C\|f\|_{W^{2,2}}^2,
\end{eqnarray}
\begin{eqnarray}
  \left[\left|f\phi_j\langle\nabla f, \nabla\mathcal{L}\phi_j\rangle\right|\right]&\leq& C\left[|x||f||\nabla f|\right]+C\left[|x|^3|f||\nabla f|\right]\nonumber\\
&\leq& C\left[|x|^2f^2\right]+C\left[|\nabla f|^2\right]+C\left[|x|^4f^2\right]+C\left[|x|^2|\nabla f|^2\right]\nonumber\\
&\leq& C\|f\|_{W^{2,2}}^2,
\end{eqnarray}
\begin{eqnarray}
  \left[\left|f\mathcal{L}\phi_j\langle\nabla f, \nabla \phi_j\rangle\right|\right]\leq C\left[|x||f||\nabla f|\right]+C\left[|x|^3|f||\nabla f|\right]\leq C\|f\|_{W^{2,2}}^2,
\end{eqnarray}
\begin{eqnarray}
  \left[f^2\left|\langle\nabla \phi_j, \nabla\mathcal{L}\phi_j\rangle\right|\right]\leq C\left[|x|^2f^2\right]+C\left[|x|^4f^2\right]\leq C\|f\|_{W^{2,2}}^2,
\end{eqnarray}
It is easy to check that
\begin{eqnarray*}
  \langle\nabla f, \nabla\langle\nabla\phi_j, \nabla f\rangle\rangle=\nabla^2\phi_j(\nabla f, \nabla f)+\nabla^2 f(\nabla f, \nabla\phi_j)
\end{eqnarray*}
and
\begin{eqnarray*}
  \langle\nabla\phi_j, \nabla\langle\nabla\phi_j, \nabla f\rangle\rangle=\nabla^2\phi_j(\nabla\phi_j, \nabla f)+\nabla^2f(\nabla\phi_j, \nabla\phi_j).
\end{eqnarray*}
It follows that
\begin{eqnarray}
  \left[\left|\phi_j\langle\nabla f, \nabla\langle\nabla\phi_j, \nabla f\rangle\rangle\right|\right]&\leq& C\left[|\phi_j||\nabla^2\phi_j||\nabla f|^2\right]+C\left[|\phi_j||\nabla^2 f||\nabla f||\nabla\phi_j|\right]\nonumber\\
&\leq&C\left[|\nabla f|^2\right]+C\left[|x|^2|\nabla f|^2\right]+C\left[|x||\nabla f||\nabla^2 f|\right]\nonumber\\
&\leq&C\left[|\nabla f|^2\right]+C\left[|x|^2|\nabla f|^2\right]+C\left[|\nabla^2 f|^2\right]\nonumber\\
&\leq&C\|f\|_{W^{2,2}}^2,
\end{eqnarray}
and
\begin{eqnarray}
  \left[\left|f\langle\nabla\phi_j, \nabla\langle\nabla\phi_j, \nabla f\rangle\rangle\right|\right]&\leq& C\left[|f||\nabla\phi_j||\nabla f||\nabla^2\phi_j|\right]+C\left[|f||\nabla^2f||\nabla\phi_j|^2\right]\nonumber\\
&\leq&C\left[|x||f||\nabla f|+|x|^3|f||\nabla f|\right]+C\left[|x|^2|f||\nabla^2f|\right]\nonumber\\
&\leq&C\left[|x|^2f^2\right]+C\left[|\nabla f|^2\right]+C\left[|x|^4f^2\right]+C\left[|x|^2|\nabla f|^2\right]\nonumber\\
&&+C\left[|x|^4f^2\right]+C\left[|\nabla^2 f|^2\right]\nonumber\\
&\leq&C\|f\|_{W^{2,2}}^2.
\end{eqnarray}
\begin{eqnarray}\label{Fj2}
  \left[\left|\langle fJ\nabla\phi_j, y\rangle\right|\right]&\leq&\left[\right(\phi_j-1)^2|\nabla f|^2]^{\frac{1}{2}}\left[|y^\perp|^2\right]^{\frac{1}{2}}+\left[f^2\left|\nabla\phi_j\right|^2\right]
  ^{\frac{1}{2}}\left[|y^\perp|^2\right]^\frac{1}{2}\nonumber\\
  &\leq& C\|f\|_{W^{1,2}}\left[|y^\perp|^2\right]^\frac{1}{2}
\end{eqnarray}
By (\ref{Fj1})-(\ref{Fj2}), and noticing the fact that $\phi_j\to 1$, $|\nabla\phi_j|\to0$, $\mathcal{L}\phi_j\to 0$, $|\nabla\mathcal{L}\phi_j|\to 0$, and $|\nabla^2\phi_j|\to 0$ pointwise as $j\to\infty$, the dominated convergence theorem gives
\begin{eqnarray*}
  \lim_{j\to\infty}\left[-\langle \nabla f_j, \nabla\left({\mathcal{L}}f_j+f_j\right)\rangle\right]&=&\left[-\langle \nabla f, \nabla\left({\mathcal{L}}f+f\right)\rangle\right]=\left[(\mu-1)|\nabla f|^2\right]\\&=&\left[\mu(\mu-1)f^2\right]<0,
\end{eqnarray*}
and for each $y$,
\begin{eqnarray*}
  \lim_{j\to\infty}\left[\langle J\nabla f_j, y\rangle\right]=0.
\end{eqnarray*}
More precisely, for any $\delta>0$, we can choose $j$ sufficiently large, such that
\begin{eqnarray*}
  \left[-\langle \nabla f_j, \nabla\left({\mathcal{L}}f_j+f_j\right)\rangle\right]\leq\frac{1}{2}\left[\mu(\mu-1)f^2\right]<0,
\end{eqnarray*}
and
\begin{eqnarray*}
  \left[\left|\langle J\nabla f_j, y\rangle\right|\right]\leq\delta\left[|y^\perp|^2\right]^\frac{1}{2}.
\end{eqnarray*}
Therefore, for any $h\in \textbf{R}$ and $y\in \textbf{R}^{2n}$, we have
\begin{eqnarray*}
  F_j''&\leq&\frac{1}{2}\mu(\mu-1)\left[f^2\right]+\delta\left[|y^\perp|^2\right]^{\frac{1}{2}}-h^2\left[|H|^2\right]-\frac{1}{2}\left[|y^\perp|^2\right]\\
  &\leq&\frac{1}{2}\mu(\mu-1)\left[f^2\right]+\frac{1}{2}\delta^2-\frac{1}{2}\left(\left[|y^\perp|^{2}\right]^\frac{1}{2}-\delta\right)^2\\
  &\leq&\frac{1}{2}\mu(\mu-1)\left[f^2\right]+\frac{1}{2}\delta^2.
\end{eqnarray*}
Choosing $\delta$ such that $\delta^2\leq\frac{1}{4}\mu(1-\mu)\left[f^2\right]$, we get $F_j''<0$ for every $h$ and $y$. This contradicts with (i).

This completes the proof of the theorem. \hfill Q.E.D.

Since the eigenvalues and eigenspaces of $\mathcal{L}$ on $\underbrace{S^1(\sqrt{2})\times\cdots S^1(\sqrt{2})}_k\times \textbf{R}^{n-k}$ ($0\leq k\leq n$) satisfy (ii), it follows immediately that
\begin{corollary}
  $\underbrace{S^1(\sqrt{2})\times\cdots S^1(\sqrt{2})}_k\times \textbf{R}^{n-k}$$\ \subset\textbf{C}^n$ ($0\leq k\leq n$) is Hamiltonian F-stable. In particular, the cylinder $S^1(\sqrt{2})\times\textbf{R}^{n-1}\subset \textbf{C}^n$ is Hamiltonian F-stable.
\end{corollary}

\begin{remark}
  In \cite{LL} (Theorem 5), Lee-Lue gave an equivalent condition for F-stability of self-shrinkers in higher codimension, which has some kind of connection with our theorem for Hamiltonian F-stability of Lagrangian self-shrinkers. In their paper, they pinned down the stability of self-shrinkers in higher codimension to the mean curvature vector being the first vector-valued eigenfunction for an elliptic system.
\end{remark}

  Note that for the special case $n=2$, by (\ref{Lx2}) we have $|x|^2-4$ is an eigenfunction of $\mathcal{L}$ which corresponds to the eigenvalue 1 if $|x|^2-4\neq0$. Besides, by the remark after Theorem 1.1 of \cite{CL}, $|x|^2=4$ implies that $\Sigma$ must be the Clifford torus which satisfies $\lambda_2(\mathcal{L})=1$. Therefore in this situation, the condition $\lambda_2({\mathcal{L}})\geq 1$ in (ii) of Theorem \ref{complete} can be changed to $\lambda_2({\mathcal{L}})=1$. Thus we get an improved statement of the characterization theorem for the $n=2$ case:

\begin{theorem} 
    Let $\Sigma^2\subset \textbf{C}^2$ be a smooth complete Lagrangian self-shrinker without boundary and with polynomial volume growth. Suppose there exist constants $C_0>0$ and $\varepsilon<\frac{1}{32}$ such that $|A|^2\leq C_0+\varepsilon |x|^2$. Then the following statements are equivalent:

  (i) $\Sigma$ is Hamiltonian F-stable.

  (ii) $\lambda_1({\mathcal{L}})=\frac{1}{2}$, $\lambda_2({\mathcal{L}})= 1$, and the eigenspace corresponding to the eigenvalue $\frac{1}{2}$ is spanned by coordinate functions.
\end{theorem}

\end{document}